\documentclass{amsart}

\usepackage{dsfont}
\usepackage{tikz-cd}
\usepackage{enumitem}
\usepackage{texmaX}

\def\triv{\mathds{1}}
\def\cP{\mathcal{P}}
\def\XXX{\mathcal{S}}

\begin{document}

\title{A note on the parity conjecture and base change}
\author{Vladimir Dokchitser}

\address{University College London, London WC1H 0AY, UK}
\email{v.dokchitser@ucl.ac.uk}


\begin{abstract}
The parity conjecture predicts that the parity of the rank of an abelian variety is determined by its global root number, that is by the sign in the conjectural functional equation of its $L$-function.
Assuming the Shafarevich--Tate conjecture, we show that if a semistable principally polarised abelian variety $A/\Q$ satisfies the parity conjecture over $\Q$ and over all quadratic fields, then it satisfies it over all number fields. More generally, we establish a criterion for when the parity conjecture for the base change of an abelian variety to a larger number field is already implied by the parity conjecture over the ground field and over other small extensions.
\end{abstract}

\maketitle

\section{Introduction}

The Birch--Swinnerton-Dyer conjecture classically predicts that the parity of the rank of the Mordell--Weil group of an abelian variety over a number field $A/K$ agrees with the order of vanishing of its $L$-function $L(A,s)$ at the central point $s\!=\!1$. In view of the conjectural functional equation $L(A,s)\leftrightarrow \pm L(A,2-s)$, the parity of the order of vanishing is determined by the sign $\pm$, which is given by the global root number $w(A/K)=\pm 1$. The resulting ``parity conjecture'' thus predicts
$$
 (-1)^{\rk A/K} = w(A/K).
$$
An advantage over the full Birch--Swinnerton-Dyer conjecture is that this formulation does not require the knowledge of the analytic continuation or the functional equation of $L(A/K,s)$, since the root number is defined independently of these. Nonetheless, almost all the results on the parity conjecture rely on some form of finiteness of Tate--Shafarevich groups, as the rank remains an impossibly difficult invariant to control otherwise.

In this note we establish the following compatibilty of the parity conjecture with extensions of the base field:

\begin{theorem}[see Theorem \ref{thm:parityextensions}, Corollaries \ref{PCcorollary2}, \ref{PCcorollary1}]
\label{thm:introextraparity}
Let $A$ be a principally polarised abelian variety over a number field $K$. Let $F/K$ be a finite Galois extension such that $\sha(A/F)$ is finite. Suppose that either
\begin{itemize}
\item $A$ is an elliptic curve and $K$ is totally real, or
\item $A$ is an elliptic curve that satisfies the parity conjecture over $K$, or
\item $A$ is semistable and satisfies the parity conjecture over $K$ and over all the quadratic extensions of $K$ in $F$.
\end{itemize}
Then the parity conjecture holds for $A$ over all extensions of $K$ in $F$.
\end{theorem}

The proof of the above result makes use of known cases of the parity conjecture for ``Artin twists'' of abelian varieties. 
The group of $F$-rational points \hbox{$A(F)\otimes_\Z\C$} is naturally a $\Gal(F/K)$-representation; we will write $A(F)_\C$ for its character. 
Assuming the finiteness of $\sha$, there is a known supply of characters $\rho$ for which one can control the parity of the inner product $\langle\rho,A(F)_\C\rangle$ and relate it to root numbers. This includes all characters of the form $\rho\!=\!\tau\! +\! \bar\tau$ and, most crucially, characters $\rho\!=\!\tau\!+\!\triv\!+\!\det\tau$ when $\tau$ is irreducible of degree 2 and $\Gal(F/K)$ is the dihedral group $D_{2p}$ for an odd prime $p$ (with suitable analogues for $C_2\times C_2$ and $D_8$).
Our contribution is the observation that this supply is sufficient to prove the above theorem. The key new input is the following purely representation-theoretic induction theorem.

\begin{theorem}\label{thm:vind}
Let $G$ be a finite group and $\cP$ a set of real-valued generalised characters of $G$ that satisfies:
\begin{enumerate}[leftmargin=*]
\item $\tau+\bar\tau\in \cP$ for every generalised character $\tau$ of degree 0, 
\item $\Ind_H^G(\tau-\triv_H-\det\tau)\in \cP$ for every subgroup $H$ and character $\tau$ of degree 2 of $H$ 
that factors through a quotient of $H$ isomorphic to $C_2\!\times\! C_2, D_8$ or to $D_{2p}$ for an odd prime~$p$,
\item if $\tau_i\in\cP$ and $n_i\in\Z$ then $\sum_i n_i \tau_i\in\cP$.
\end{enumerate}
Then $\cP$ contains all permutation characters of $G$ that have degree 0 and trivial determinant. 
\end{theorem}

Here by a {\em generalised character} of a finite group $G$ we mean a $\Z$-linear combination of irreducible complex characters of $G$.
A {\em permutation character} is a $\Z$-linear combination of characters of permutation representations.
We extend the notion of the {\em determinant} of a character of a representation to all characters by setting $\det (\tau_1-\tau_2)=\det\tau_1 \otimes \det\bar{\tau}_2$.
$D_{2n}$ denotes the dihedral group of order $2n$.

In the context of the parity conjecture, it is crucial that the induction theorem only makes use of the most basic dihedral groups. Indeed, to the best of the author's knowledge, there is yet no way to control the parity of the inner product $\langle \tau\!+\!\triv\!+\!\det\tau,A(F)_\C\rangle$ for general dihedral groups: neither the techniques of \cite{squarity, tamroot} nor of \cite{MazurRubin} appear to have been extended to these cases. (The most elementary problematic case is $D_{42}$.)

\begin{acknowledgements}
The author would like to thank Holly Green, Alexandros Konstantinou and Adam Morgan for many useful discussions and for our joint project \cite{brauer}. Indeed, the main motivation for the above induction theorem came from that work, where it is applied to Galois groups of covers of curves in order to control the parity of ranks of Jacobians.
\end{acknowledgements}

\section{Proof of the induction theorem}

In this section we prove Theorem \ref{thm:vind}.

\begin{notation}
We write $\XXX_G$ for the set of permutation characters of $G$ that have degree 0 and trivial determinant. 

We write $\cP_G$ for the smallest set of generalised characters closed under (1), (2), (3) of Theorem \ref{thm:vind}. Equivalently, $\rho\in\cP_G$ if and only if $\rho=\sum_i n_i \rho_i$ with each $\rho_i$ either of the form 
\begin{itemize}
  \item $\rho_i=\tau_i+\bar\tau_i$ for a generalised character $\tau$ of degree 0, or 
  \item $\rho_i=\Ind_{H_i}^G(\tau_i-\triv_{H_i}-\det\tau_i)$ for some subgroup $H_i$ and degree 2 character $\tau_i$ of $H_i$ that factors through a quotient of $H_i$ isomorphic to $C_2\!\times\! C_2, D_8$ or $D_{2p}$ for an odd prime~$p$.
\end{itemize}
\end{notation}

\begin{lemma}\label{lem:indx}
Let $G$ be a finite group.
\begin{enumerate}
\item If $H\le G$ and $\rho\in \XXX_H$ then $\Ind_H^G\rho\in \XXX_G$.
\item If $H\le G$ and $\rho\in \cP_H$ then $\Ind_H^G\rho\in \cP_G$.
\item If $N\triangleleft G$ then $\XXX_{G/N}\subseteq \XXX_G$ and $\cP_{G/N}\subseteq \cP_G$.
\item All $\rho\in\cP_G$ have degree 0 and trivial determinant.
\end{enumerate}
\end{lemma}

\begin{proof}
(1) $\Ind_H^G\rho$ is clearly a permutation character of degree 0.

It thus suffices to check that if $\rho$ has degree 0 and trivial determinant, then $\Ind_H^G\rho$ has trivial determinant. From the explicit construction of induced representations, it follows that for a character $\tau$ of $H$
$$
\det (\Ind_H^G\tau)(g) = \epsilon^{\deg\tau}\prod_{i=1}^{[G:H]}\det \tau(h_i),
$$
for some suitable $h_i\in H$ and $\epsilon\in \pm 1$ (the sign of the permutation of $g$ on $G/H$). Thus if $\det\tau=\triv$ and the degree of $\tau$ is even, then $\det\Ind\tau$ is also trivial. This readily extends to generalised characters.

(2) Clear.

(3) Clear.

(4) The degree 0 claim is clear. The determinant claim follows as in the proof of~(1).
\end{proof}

\begin{lemma}\label{lem:abeliannice}
If $\rho\in \XXX_G$ is a $\Z$-linear combination of linear characters, then $\rho\in\cP_G$.
\end{lemma}

\begin{proof}
The only real-valued linear characters are either trivial or have order 2, so we can write $\rho = (\tau+\bar\tau)+(\sum_{i=1}^k \epsilon_i)+m\triv$ for some generalised character $\tau$, some linear characters $\epsilon_i$ of order 2, and $m=-2\deg\tau-k$. 

For any two linear order 2 characters $\epsilon, \epsilon'$, we have $\epsilon+\epsilon'-\epsilon\otimes\epsilon'-\triv \in \cP_G$, as it is either 0 or in $\cP_{C_2\times C_2}$ for a $C_2\times C_2$-quotient of $G$. Hence $\sum\epsilon_i - \otimes_i\epsilon_i-(k-1)\triv\in \cP_G$.
Also, $\tau+\bar\tau - (2\deg\tau) \triv \in \cP_G$.
We can therefore write $\rho=\rho'+\bigotimes_i\epsilon_i-n\triv$ for some $\rho'\in \cP_G$ and $n\in\Z$. Since $\det\rho=\det\rho'=\triv$ and $\rho$ and $\rho'$ have degree 0 by Lemma \ref{lem:indx}(4), we must have $n=1$ and $\bigotimes_i\epsilon_i=\triv$. Hence $\rho=\rho'\in\cP_G$.
\end{proof}

\begin{lemma}\label{lem:oddnice}
If $G$ has odd order then $\XXX_G\subseteq \cP_G$.
\end{lemma}

\begin{proof}
In a group of odd order, the only real-valued irreducible character is $\triv$. So every real-valued generalised character is of the form $\tau+\bar\tau+m\triv$ for some $m\in\Z$.
\end{proof}

\begin{lemma}\label{lem:rhoH}
Every $\rho\in \XXX_G$ can be written in the form 
$$
\rho=\sum_H \rho_H+ (\text{a sum of linear characters})
$$ for some set of subgroups $H\le G$, where
$$
  \rho_H=\Ind_H^G\triv - \det\Ind_H^G\triv - ([G\!:\!H]-\!1)\triv \>\>\in \XXX_G.
$$ 
\end{lemma}
\begin{proof}
By construction, $\rho_H$ is a permutation character of degree 0 and with trivial determinant, and is of the form $\Ind_H^G\triv_H + $(a sum of linear characters).
Permutation characters are linear combinations of $\Ind_H^G\triv$, so subtracting the corresponding $\rho_H$'s from $\rho$ leaves a sum of linear characters.
\end{proof}

\pagebreak

\begin{proposition}\label{prop:sylownice}
If the $p$-Sylow subgroup $N$ of $G$ is normal and cyclic, and $\XXX_{G/N}\subseteq \cP_{G/N}$, then $\XXX_G\subseteq \cP_G$.
\end{proposition}

\begin{proof}
\underline{Case 1: $p\!=\!2$}. 
By the Schur--Zassenhaus theorem, we can write $G=N\rtimes G_0$ with $|G_0|$ odd. Since $N$ is cyclic of 2-power order, $\Aut N$ is a 2-group, and hence $G=N\times G_0$. The only real-valued irreducible characters of this group are $\triv$ and $\epsilon$, the order 2 linear character that factors through $G/G_0\simeq N$. Thus $\rho\in \XXX_G$ is of the form $\rho=\tau+\bar{\tau}+a\epsilon+b\triv$. The determinant and degree conditions ensure that $a$ and $b$ are even, so $\rho= (\tau+\frac{a}{2}\epsilon+\frac{b}{2}\triv)+\overline{(\tau+\frac{a}{2}\epsilon+\frac{b}{2}\triv)}\in\cP_G$.

\underline{Case 2: $p$ is odd}. 
We proceed by induction on $|G|$.

By Lemmata \ref{lem:abeliannice} and \ref{lem:rhoH} it suffices to show that $\rho_H\in\cP_G$ for all $H\le G$, where $\rho_H\in \XXX_G$ is given by
$$
  \rho_H=\Ind_H^G\triv_H - \det\Ind_H^G\triv_H - ([G\!:\!H]-\!1)\triv_G.
$$
If $N$ is trivial then there is nothing to prove. So suppose that $N$ is non-trivial and let $V=C_p\le N$. As $V$ is characteristic in $N$ and $N$ is normal in $G$, $V$ is normal in $G$. 

Suppose $H\supseteq V$. Then, by induction, $\rho_H\in \cP_{G/V}\subseteq \cP_G$. 

Suppose $H\supsetneq V$ and $VH\neq G$. 
Working in $VH$, by induction,
$$
\rho_0=\Ind_H^{VH}\triv_H - \det\Ind_H^{VH}\triv_H - ([VH\!:\!H]-\!1)\triv_{VH} \qquad \in \cP_{VH}.
$$ 
As $VH=V\rtimes H$ and $V$ has odd order, $\det\Ind_H^{VH}\triv$ is either trivial or can be written as $\Ind_{VH'}^{VH}\triv_{VH'} - \triv_{VH}$ for some index 2 subgroup $H'$ of $H$. Thus
$$
  \Ind_{VH}^G\rho_0 = \Ind_H^G\triv_H+ a \Ind_{VH'}^{G}\triv_{VH'} + b \Ind_{VH}^{G}\triv_{VH} 
$$
for some suitable $a, b\in\Z$.
We can therefore write
$$
\rho_H = \Ind_{VH}^G\rho_0+a\rho_{VH'}+b\rho_{VH}+{\text{(sum of linear characters)}}.
$$
By Lemma \ref{lem:indx}, $\Ind\rho_0\in\cP_G$. Moreover $\rho_{VH'}, \rho_{VH} \in \cP_G$ by the ``$H\supseteq V$'' case. Finally, $\rho_H, \rho_{VH}, \rho_{VH'}\in \XXX_G$ by construction and $\Ind\rho_0\in \XXX_G$ by Lemma \ref{lem:indx}, so Lemma \ref{lem:abeliannice} shows that $\rho_H\in\cP_G$.

Suppose that $H\supsetneq V$ and $VH=G$ and $W=\ker (H\to \Aut V)$ is non-trivial. In this setting $G=V\rtimes H$ and $W\triangleleft G$. So, by induction and Lemma \ref{lem:indx}, $\rho_H\in\cP_{G/W}\subseteq \cP_G$.

Finally, suppose that $H\supsetneq V$, $VH=G$ and $H$ acts faithfully on $V$ by conjugation. 
In other words,
$G=\F_p\rtimes H$ with $H\le \F_p^\times$. If $H$ has odd order, the result follows from Lemma~\ref{lem:oddnice}.
Otherwise $C_2\le H$ and $\F_p\rtimes C_2=D_{2p}\lhd G$. 
The irreducible characters of this group are either linear ones that factor through $G/\F_p$, or those of the form $\Ind_{D_{2p}}^G\tau$ for a 2-dimensional irreducible character $\tau$ of $D_{2p}$. 
Noting that $\Ind_{D_{2p}}^G\det\tau$ and $\Ind_{D_{2p}}^G\triv$ are sums of linear characters of $G$,
we obtain
$$
\rho_H = \sum_i \Ind_{D_{2p}}^G(\tau_i-\det\tau_i-\triv_{D_{2p}}) + (\text{sum of linear characters})
$$
for some 2-dimensional irreducible characters $\tau_i$ of $D_{2p}$.
As $\rho_H\in \XXX_G$ and $\sum_i \Ind_{D_{2p}}^G(\tau_i-\det\tau_i-\triv_{D_{2p}})\in \cP_G$, they both have degree 0, trivial determinant and real character.
Hence so does the final sum of linear characters of the cyclic group $G/\F_p\simeq H$, which must therefore be of the form $\sigma+\bar\sigma\in \cP_G$. Hence $\rho_H\in\cP_G$ as well.
\end{proof}

\begin{lemma}\label{lem:solomon}
For all $\rho\in \XXX_G$, there are 
hyperelementary\footnote{Recall that a hyperelementary group is one of the form $C_n\rtimes P$, where $P$ is a $p$-group} 
$H_i\le G$, $\rho_i\in \XXX_{H_i}$ and $n_i\in\Z$ such that
$$
\rho = \sum_i n_i\Ind_{H_i}^G\rho_i
$$
\end{lemma}

\begin{proof}
By Solomon's Induction theorem (\cite{Solomon} Thm.\ 1), there are hyperelementary $H_i\le G$ and $n_i\in \Z$ with $\triv_G=\sum_i\Ind_{H_i}^G\triv_{H_i}$. Thus
$$
\rho = \sum_i n_i (\rho\otimes \Ind_{H_i}^G\triv_{H_i})=\sum_i n_i \Ind_{H_i}^G \Res_{H_i}^G \rho.
$$
The result follows as $\rho_i=\Res_{H_i}^G \rho\in \XXX_{H_i}$.
\end{proof}

\begin{theorem}
For every finite group $\XXX_G\subseteq \cP_G$.
\end{theorem}

\begin{proof}
\underline{Case 1: $G$ has odd order}. 
Lemma \ref{lem:oddnice}.

\underline{Case 2: $G$ is a 2-group}.
Let $\rho\in \XXX_G$.

By Lemmata \ref{lem:abeliannice} and \ref{lem:rhoH} we may assume that 
$$
\rho=\Ind_H^G\triv + (\text{a sum of linear characters}),
$$
for some $H\le G$.

We proceed by induction on $[G\!:\!H]$. When $[G\!:\!H]\le 2$ the result follows from Lemma~\ref{lem:abeliannice}.
If $[G\!:\!H]\ge 4$, pick subgroups $H\le U \le V\le G$ with $[U\!:\!H]=[V\!:\!U]=2$, which exist as $G$ is a 2-group.
If $H\lhd V$ and $V/H\simeq C_4$, then $\Ind_H^G\triv=\Ind_U^G\triv+(\tau+\bar\tau)$ for some $\tau$, and the result follows by induction.
If $H\lhd V$ and $V/H\simeq C_2\times C_2$ and $U, U', U''$ are the three intermediate subgroups, then $\Ind_H^G\triv=\Ind_U^G\triv+\Ind_{U'}^G\triv+\Ind_{U''}^G\triv-2\Ind_V^G\triv$, and the result follows by induction. Finally, if $H$ is not normal in $V$ then there is an index 2 subgroup $H_0\le H$ with $H_0\lhd V$ with $V/H_0\simeq D_8$; in this case, writing $\sigma$ for the 2-dimensional irreducible representation of $D_8$, 
$$
\Ind_H^G\triv=\Ind_{H_0}^G(\sigma+\triv+\det\sigma)=\Ind_{H_0}^G(\sigma-\triv-\det\sigma)+(\tau+\bar\tau)+4[G\!:\!H_0]\triv_G
$$ 
for $\tau=\Ind_{H_0}^G(\triv+\det\sigma)-2[G\!:\!H_0]\triv_G$, so that $\rho\in \cP_G$ by Lemma~\ref{lem:abeliannice}.

\underline{Case 3: $G$ is hyperelementary}.
Write $G=C_n\rtimes P$, where $P$ is a $p$-group and $n$ is coprime to $p$. By Cases 1 and 2, $\XXX_{G/C_n}\subseteq \cP_{G/C_n}$. The result follows by a repeated application of Proposition \ref{prop:sylownice}.

\underline{Case 4: $G$ arbitrary}.
Combine Lemma \ref{lem:solomon}, Lemma \ref{lem:indx} and Case 3.
\end{proof}

This completes the proof of Theorem \ref{thm:vind}.
We end this section by recording an alternative version of the theorem.

\begin{corollary}
Let $G$ be a finite group and $\cP$ a set of real-valued characters of $G$ that satisfies:
\begin{itemize}[leftmargin=*]
\item $\Ind_H^G\tau\in \cP$ for every subgroup $H$ and real character $\tau$ of degree 0 and trivial determinant of $H$ that factors through a quotient of $H$ that is either cyclic or isomorphic to $C_2\!\times\! C_2, D_8$ or to $D_{2p}$ for some odd prime~$p$,
\item if $\tau_i\in\cP$ and $n_i\in\Z$ then $\sum_i n_i \tau_i\in\cP$.
\end{itemize}
Then $\cP$ contains all permutation characters of $G$ that have degree 0 and trivial determinant. 
\end{corollary}

\begin{proof}
$\cP$ clearly satisfies (2) and (3) of Theorem \ref{thm:vind}.

By the degree 0 variant of Brauer's induction theorem (see \cite{DelC} Prop.\ 1.5), a character $\tau$ of degree 0 of $G$ can be written as 
$$
\tau=\sum_i \Ind_{H_i}^{G} \psi_i
$$ 
for some subgroups $H_i\le G$ and degree 0 characters $\psi_i$ of $H_i$ that are $\Z$-linear combinations of linear characters.  
Explicitly writing out $\psi_i = \sum_j n_{i,j} \chi_{i,j}$ for linear characters $\chi_{i,j}$, and using the fact that $\sum_j n_{i,j}=\deg\psi_i=0$, we get
$$
 \tau+\bar\tau = \sum_i \Ind_{H_i}^{G} (\psi_i + \bar\psi_i) = 
\sum_i \Ind_{H_i}^{G}  \sum_j n_{i,j} (\chi_{i,j} + \bar\chi_{i,j} - 2\triv).
$$
Since $\chi_{i,j} + \bar\chi_{i,j} - 2\triv$ is real, has degree~0, trivial determinant and factors through a cyclic quotient of $H_i$, this shows that $\cP$ also satisfies (1)  of Theorem \ref{thm:vind}.
\end{proof}

\section{Applications to the parity conjecture}

We now turn to the behaviour of the parity conjecture under base change.

Let $A/K$ be an abelian variety over a number field, and $F/K$ a finite Galois extension. To a real character of $\tau$ of $\Gal(F/K)$ one associates a root number $w(A/K, \tau)\in \{\pm 1\}$, which , conjecturally, determines the sign in the functional equation for the twisted $L$-function $L(A/K,s)$. One also expects the following ``parity conjecture for twists'' (see e.g. \cite{tamroot} \S1.1 or the introduction to \cite{RohG}):
$$
 (-1)^{\langle \tau, A(F)_\C \rangle} = w(A/K,\tau).
$$
As before, $A(F)_\C$ is the character of $A(F)_\Z\otimes\C$, the $\Gal(F/K)$-representation on the $F$-rational points of $A$. (Here we use the natural extension of the definition of root numbers of twists from characters of representations to generalised characters by $w(A/K,\tau_1-\tau_2)=\frac{w(A/K,\tau_1)}{w(A/K,\tau_2)}$.)

The parity conjecture for twists is known in the following cases for dihedral Galois groups.

\begin{theorem}\label{thm:knowntwists}
Let $A$ be a principally polarised abelian variety over a number field $K$ and $p$ a prime number.
Let $F/K$ be a finite Galois extension with Galois group $\Gal(F/K)\simeq C_2\!\times\! C_2$ or $D_8$ if $p=2$, or $\Gal(F/K)\simeq D_{2p}$ if $p\neq 2$.
If $A/K$ is either an elliptic curve or it has semistable reduction at all primes that ramify in $F/K$, and if the $p$-primary component of $\sha(A/F)$ is finite, then for every character $\tau$ of $\Gal(F/K)$ of degree 2
$$
 (-1)^{\langle \tau -\triv-\det\tau, \>A(F)_\C \rangle}= w(A/K,\tau -\triv-\det\tau).
$$
In other words, the parity conjecture holds for the twist of $A/K$ by  $\tau -\triv-\det\tau$ .
\end{theorem}

\begin{proof}
For $A/K$ satisfying the semistability hypothesis this is \cite{tamroot} Thm.\ 4.2 ($p\neq 2$) and \cite{betts} Thm.\ 1.3.2 applied to \cite{tamroot} Ex.\ 2.54 ($p=2$).
For elliptic curves this is \cite{kurast} Thm.\ 6.7 ($p=2, 3$) and \cite{Rochefoucauld} Thm.\ 2.1 ($p\ge 5$).
\end{proof}

We now have the necessary tools for the proof of Theorem \ref{thm:introextraparity}.
In what follows, we will freely use the standard fact that for an abelian variety over a number field $A/K$ and a finite extension $F/K$, if $\sha(A/F)[p^\infty]$ is finite then so is $\sha(A/K)[p^\infty]$; see e.g. \cite{squarity} Remark 2.10.

\begin{theorem}\label{thm:parityextensions}
Let $A$ be a principally polarised abelian variety over a number field~$K$. 
Let $F/K$ be a finite Galois extension such that $\sha(A/F)[p^\infty]$ is finite for all $p\,|\,[F\!:\!K]$. 
Suppose moreover that either $A/K$ has semistable reduction at all primes that ramify in $F/K$, or that $A$ is an elliptic curve. 
If the parity conjecture holds for $A$ over $K$ and over all the quadratic extensions of $K$ in $F$, then it holds over all intermediate fields of $F/K$. 
\end{theorem}

\begin{proof}
First recall the following basic properties of root numbers of Artin twists of abelian varieties (the first and the final one come direct from the definitions, for the others see e.g.  \cite{tamroot} Proposition A.2):
$$
\begin{array}{rcl}
w(A/K,\rho_1+\rho_2)&=& w(A/K,\rho_1)w(A/K,\rho_2), \cr
w(A/K,\bar{\rho}+\rho)&=&1,\cr
w(A/K,\Ind_H^G\rho)&=&w(A/F^H,\rho),\cr
w(A/K,\triv) &=& w(A/K),
\end{array}
$$
where $G=\Gal(F/K)$ and $H\le G$ is a subgroup.

Write $A(F)_\C$ for the character of $A(F)\otimes_\Z\C$ as a $G$-representation. 

The character $A(F)_\C$ is clearly $\Q$-valued, so for  $\rho=\tau+\bar{\tau}$ we have $(-1)^{\langle \rho, A(F)_\C\rangle}=1=w(A,\rho)$.
If $\rho=\Ind_H^G(\tau-\triv_H-\det\tau)$ for a character $\tau$ of $H$ with $\deg\tau=2$ that factors through a quotient of $H$ isomorphic to $C_2\!\times\! C_2, D_8$ or $D_{2p}$ for an odd prime~$p$, then 
$$
(-1)^{\langle \rho, A(F)_\C\rangle}=(-1)^{\langle \tau-\triv_H-\det\tau,  \>A(F)_\C\rangle_H} = w(A/F^H,\tau-\triv_H-\det\tau)= w(A/K,\rho)
$$ 
by Frobenius reciprocity and Theorem \ref{thm:knowntwists}.
Hence, by Theorem \ref{thm:vind}, $(-1)^{\langle A(F)_\C,\rho\rangle}=w(A,\rho)$ for every permutation character $\rho$ of degree 0 and with trivial determinant.

By hypothesis and the above basic properties of root numbers, the formula  $(-1)^{\langle A(F)_\C,\rho\rangle}=w(A,\rho)$ also holds when $\rho$ is either the trivial character or a character of order 2. It follows that it holds for by all permutation characters $\rho$.

Thus, by Frobenius reciprocity, for any subgroup $H\le G$ ,
$$
(-1)^{\rk A/F^H}=(-1)^{\langle \Ind_H^G\triv, A(F)_\C \rangle}=w(A,\Ind_H^G\triv)=w(A/F^H).
$$
\end{proof}

\begin{corollary}\label{PCcorollary2}
Let $E/K$ be an elliptic curve over a number field satisfying the parity conjecture, and $F/K$ a finite Galois extension. 
If $\sha(E/F)[p^\infty]$ is finite for all primes $p\, |\,[F\!:\!K]$, then the parity conjecture holds for $E$ over all intermediate fields of $F/K$.
\end{corollary}

\begin{proof}
The parity conjecture holds for $E$ over every quadratic extension of $K$ inside $F$ by the Kramer--Tunnell theorem (\cite{KT} and \cite{kurast} Corollary 1.6) and the assumed finiteness of the 2-primary part of $\sha$ over that extension. 
\end{proof}

\begin{corollary}\label{PCcorollary1}
Let $E/K$ be an elliptic curve over a totally real field, and let $F/K$ be a finite Galois extension. 
If $\sha(E/F)[p^\infty]$ is finite for all primes $p\, |\,[F\!:\!K]$, then the parity conjecture holds for $E$ over all intermediate fields of $F/K$. (If $F\!=\! K$ we require finiteness for some $p$.)
\end{corollary}

\begin{proof}
The parity conjecture holds for $E/K$ 
by the $p$-parity conjecture for elliptic curves over totally real fields (\cite{Nekovar3} Thm.\ E and \cite{HollyCeline} Thm.\ 1.5) and the assumed finiteness of $\sha(E/K)[p^\infty]$ for some $p$.
\end{proof}

\begin{corollary}\label{PCcorollary3} 
Let $K$ be a number field and let $J/K$ be the Jacobian of a semistable hyperelliptic curve of genus $g\geq 2$ that has good reduction at all primes above 2 and that satisfies the parity conjecture over $K$. 
Let $F/K$ be a Galois extension in which all primes of $K$ above 2 are unramified.
If $\sha(J/F)[p^\infty]$ is finite for all primes $p\, |\,[F\!:\!K]$, then the parity conjecture holds for $J$ over all intermediate fields~of~$F/K$.
\end{corollary}

\begin{proof}
The parity conjecture holds for $J$ over a quadratic extension of $K$ inside $F$ by \cite{morgan} Thm.\ 1.1
and the assumed finiteness of the 2-primary part of $\sha$ over that extension. 
\end{proof}

\begin{remark}
The assumption on the finiteness of $\sha(A/F)[p^\infty]$ in Theorem \ref{thm:parityextensions} can be weakened from ``finiteness for all $p\,|\,[F\!:\!K]$'' to ``finiteness for all odd $p$ for which $\Gal(F/K)$ contains a $D_{2p}$-subquotient and for $p=2$ if $\Gal(F/K)$ contains a $C_2\!\times\! C_2$-subquotient''. This follows directly from the proof. 
\end{remark}



\end{document}